\documentclass[11pt]{amsart}

\usepackage{amssymb}  

\theoremstyle{plain}
\newtheorem{theorem}{Theorem}[section]
\newtheorem{lemma}[theorem]{Lemma}
\newtheorem{proposition}[theorem]{Proposition}
\newtheorem{corollary}[theorem]{Corollary}
\newtheorem{problem}[theorem]{Problem}
\theoremstyle{definition}
\newtheorem{definition}[theorem]{Definition}
\theoremstyle{remark}
\newtheorem{remark}[theorem]{Remark}
\newtheorem{remarks}[theorem]{Remarks}
\newtheorem{example}[theorem]{Example}
\newtheorem*{acknowledgment}{Acknowledgment}

\numberwithin{equation}{section}

\newcommand{\seclabel}[1]{\label{sec:#1}} 
\newcommand{\thmlabel}[1]{\label{thm:#1}} 
\newcommand{\lemlabel}[1]{\label{lem:#1}} 
\newcommand{\corolabel}[1]{\label{coro:#1}} 
\newcommand{\proplabel}[1]{\label{prop:#1}} 
\newcommand{\deflabel}[1]{\label{def:#1}} 
\newcommand{\remlabel}[1]{\label{rem:#1}} 
\newcommand{\exlabel}[1]{\label{ex:#1}} 
\newcommand{\eqlabel}[1]{\label{eq:#1}} 
\newcommand{\problabel}[1]{\label{prob:#1}} 

\newcommand{\secref}[1]{\ref{sec:#1}} 
\newcommand{\thmref}[1]{\ref{thm:#1}} 
\newcommand{\lemref}[1]{\ref{lem:#1}} 
\newcommand{\cororef}[1]{\ref{coro:#1}} 
\newcommand{\propref}[1]{\ref{prop:#1}} 
\newcommand{\defref}[1]{\ref{def:#1}} 
\newcommand{\remref}[1]{\ref{rem:#1}} 
\newcommand{\exref}[1]{\ref{ex:#1}} 
\renewcommand{\eqref}[1]{\ref{eq:#1}} 
\newcommand{\peqref}[1]{(\eqref{#1})} 

\newcommand{\LL}{\mathcal{L}}
\newcommand{\KK}{\mathcal{K}}
\newcommand{\MM}{\mathcal{M}}
\newcommand{\Mlt}{\mathrm{Mlt}}
\newcommand{\LMlt}{\mathrm{LMlt}}
\newcommand{\PMlt}{\mathrm{PMlt}}
\newcommand{\Nuc}{\mathrm{Nuc}}
\newcommand{\Ker}{\mathrm{ker}}
\newcommand{\Btp}{\mathrm{Btp}}
\newcommand{\Aut}{\mathrm{Aut}}
\newcommand{\iv}{^{-1}}

\title[On Twisted Subgroups and Bol Loops of Odd Order]
{On Twisted Subgroups and Bol Loops of Odd Order}

\author[T.~Foguel]{Tuval~Foguel}
\address{Department of Mathematics \\
Auburn University Montgomery \\
PO Box 244023\\
Montgomery, AL 36124-4023 USA}
\email{tfoguel@mail.aum.edu}
\author[M.~K.~Kinyon]{Michael~K.~Kinyon}
\address{Department of Mathematical Sciences \\
Indiana University South Bend \\
South Bend, IN 46634 USA}
\email{mkinyon@iusb.edu}
\urladdr{http://mypage.iusb.edu/\symbol{126}mkinyon}
\author[J.~D.~Phillips]{J.~D.~Phillips}
\address{Department of Mathematics and Computer Science\\
Wabash College\\
Crawfordsville, IN  47933}
\email{phillipj@wabash.edu}
\urladdr{http://www.wabash.edu/depart/math/faculty.html{\#}Phillips}


\begin{document}

\begin{abstract}
In the spirit of Glauberman's fundamental work in B-loops and Moufang
loops \cite{Gl1, Gl2}, we prove Cauchy and strong Lagrange theorems
for Bol loops of odd order. We also establish necessary conditions for the
existence of a simple Bol loop of odd order, conditions which should be useful in the
development of a Feit-Thompson theorem for Bol loops. Bol loops are closely related to
Aschbacher's twisted subgroups \cite{Asch}, and we survey the latter in some
detail, especially with regard to the so-called Aschbacher radical.
\end{abstract}

\maketitle


\section{Introduction}
\seclabel{introduction}

A \emph{magma} $(\LL,\cdot)$ consists of a set $\LL$ together with a binary operation $\cdot$
on $\LL$. For $x\in\LL$, define the left (resp., right) translation by $x$ by $L(x)y = x\cdot y$ (resp.,
$R(x)y = y\cdot x$) for all $y\in\LL$. A magma with a two-sided neutral element $1$ such that
all left translations bijective is called a \emph{left loop}. A left loop in which all right translations
are bijective is called a \emph{loop}. For basic facts about loops, we refer the reader to
\cite{Bel2, Br, CPS, Pf}. A loop satisfying the \textit{left Bol identity}
\[
(x\cdot (y\cdot x))\cdot z = x\cdot (y\cdot (x\cdot z))
\]
or equivalently
\[
L(x\cdot (y\cdot x)) = L(x)L(y)L(x) 
\]
for all $x,y,z\in\LL$, is called a \emph{left Bol loop}. A loop satisfying the mirror
identity $((x\cdot y)\cdot z)\cdot y = x\cdot ((y\cdot z)\cdot y)$ for all $x,y,z\in\LL$
is called a \emph{right Bol loop}, and a loop which is both left and right Bol is
a \emph{Moufang loop}. For the balance of this paper, the term ``Bol loop" will
refer to left Bol loop; all statements about left Bol loops dualize trivially to right
Bol loops. For basic facts about Bol loops, we refer the reader to \cite{Ro} and IV.6 in
\cite{Pf}. (In both cases translating from right Bol to left Bol). A \textit{Bruck loop}
is a Bol loop with the \textit{automorphic inverse property}, i.e.,
$x\iv\cdot y\iv = (x\cdot y)\iv$. (These are also known as K-loops
\cite{Kie} and gyrocommutative gyrogroups \cite{Ung}.) A loop is said to be
\emph{uniquely} $2$-\emph{divisible} if the squaring map $x\mapsto x\cdot x$
is a bijection; we will abuse terminology a bit and drop the ``uniquely". A $2$-divisible
Bruck loop is called a \textit{B-loop} \cite{Gl1}. (Glauberman's original definition
was restricted to the finite case.)

In the fundamental papers \cite{Gl1, Gl2}, Glauberman studied finite B-loops and
finite Moufang loops of odd order. In \cite{Gl1}, he proved Hall, Sylow, Cauchy and Lagrange
theorems for finite B-loops. In \cite{Gl2}, he used the B-loop results to establish similar results
for Moufang loops. He also proved Feit-Thompson theorems for both finite B-loops and finite
Moufang loops of odd order. This naturally raises the question as to how far these results extend to the general case
of finite Bol loops of odd order. In this paper, we begin examining this question. We
make use of the notion of a \emph{twisted subgroup} of a group, adopting the
terminology of Aschbacher \cite{Asch}. This same idea can be found in Glauberman's
papers \cite{Gl1, Gl2}, and we use his results to establish Cauchy and strong
Lagrange theorems for Bol loops of odd order. We also start an attack on a
Feit-Thompson theorem for Bol loops of odd order. We were not able to prove
a complete Feit-Thompson result, but we present some conditions that a simple Bol
loop of odd order must satisfy which we think will be crucial in a proof, if there is
indeed such a theorem. We also observe that certain varieties of Bol loops of odd
order, such as those in which every left inner mapping is an automorphism, are
necessarily solvable.

In the next section, we present a few preliminaries from loop theory.
This can be safely skipped by those who are more interested in groups
than in loops. Such readers will find \S\secref{ts} and \S\secref{2dts}
to their taste. A minimal amount of loop theory is present in
\S\secref{2dts}, so as not to abandon completely the spirit of
\cite{Gl1}, although it is possible in principle to avoid loops completely.
In \S\secref{bol} and \S\secref{2dbol}, we apply the results of \S\secref{ts}
and \S\secref{2dts}, respectively, to Bol loops of odd order.

Throughout this paper, we state several open problems in the hope of stimulating
research into Bol loops of odd order. In the Feit-Thompson direction, the existence
of \emph{any} finite simple (non-Moufang) Bol loop is widely considered to be
the most important open problem in loop theory \cite{OpenProblems}, and we
think that focusing on Bol loops of odd order is a reasonable place to start.


\section{Preliminaries}
\seclabel{preliminaries}

In this section, we review a few necessary notions from loop theory, and establish some
notation conventions. Binary operations in left loops will be explicitly denoted, while group
operations in groups will be denoted by juxtaposition. Permutations, such as
left and right translations, will act on the left of their arguments.

For a set $S$, we let $S!$ denote the group of all permutations of $S$. The \emph{multiplication group},
$\Mlt(\LL)$, of a loop $\LL$ is the subgroup of $\LL !$ generated by all right and left translations. The
\emph{left multiplication group}, $\LMlt(\LL)$, of a left loop $\LL$ is the subgroup of $\Mlt(\LL)$ generated
by left translations. The subgroup $\LMlt_1(\LL) = \{ \phi\in\LMlt(\LL) : \phi 1 = 1 \}$ is called the
\emph{left inner mapping group} of $\LL$. This subgroup has trivial core (recall that the core
$\mathrm{ker}_H(G) = \bigcap_{g\in G} gHg^{-1}$ of a subgroup $H$ in a group $G$ is the
largest normal subgroup of $G$ contained in $H$). The set $L(\LL) = \{ L(x) : x\in\LL \}$ of left
translations is a left transversal (complete set of coset  representatives) to each conjugate of
$\LMlt_1(\LL)$ in $\LMlt(\LL)$.

These observations lead us to the following construction (\cite{Baer}; see also \cite{KJ}). 
Let $G$ be a group, $H\leq G$, and $T\subseteq G$ a left transversal of $H$.
There is a natural $G$-action on $T$, which we denote by $\cdot$, defined by the equation
$(g\cdot x)H = gxH$, that is, $g\cdot x$ is the unique representative in $T$ of the coset $gxH$.
This action restricted to $T$ itself endows $T$ with a binary operation. If $1\in T$, then
$(T, \cdot)$ turns out to be a left loop, which we call the \textit{induced} left loop. If
$T$ is also a left transversal of each conjugate $gHg\iv$, $g\in G$, then $(T,\cdot)$ is
a loop. All of the induced left loops we discuss in this paper turn out to be loops.

\begin{proposition} \cite{Ph}
\proplabel{core}
Let $G$ be a group with subgroup $H\leq G$ and a left transversal $T\subset G$ of $H$
such that $\langle T\rangle = G$. The permutation representation $G\to T!$ defined by
$(g\cdot x)H = gxH$ ($g\in G$, $x\in T$) gives an epimorphism from $G$ onto $\LMlt(T,\cdot)$. 
The sequence $1 \to \mathrm{ker}_H(G) \to G \to \LMlt(T,\cdot) \to 1$ is exact.
\end{proposition}

If $\LL$ is a Bol loop, then $\LL$ is \emph{power-associative}, that is,
if $x^0 := 1$, $x^{n+1} := x\cdot x^n$, $x^{-n-1} := x\iv\cdot x^{-n}$, $n\geq 0$,
then $x^m \cdot x^n = x^{m+n}$ for all $x\in \LL$ and all integers $m,n$. Moreover,
$\LL$ is \emph{left power-alternative}, which means that $L(x^n) = L(x)^n$ for all
$x\in \LL$ and all integers $n$. Taking $n = -1$ and $n = 2$, we obtain, respectively, the
\emph{left inverse property} (LIP) $L(x)\iv = L(x\iv)$ and the \emph{left alternative
property} (LAP) $L(x)^2 = L(x^2)$.

The \emph{left nucleus}, \emph{middle nucleus}, \emph{right nucleus}, and \emph{nucleus}
of a loop $\LL$ are defined, respectively, by
\[
\begin{array}{rcl}
\Nuc_l(\LL) &:=& \{x \in \LL : x(yz) = (xy)z\ \forall y,z \in \LL \} \\
\Nuc_m(\LL) &:=& \{y \in \LL : x(yz) = (xy)z\ \forall x,z \in \LL \} \\
\Nuc_r(\LL) &:=& \{z \in \LL : x(yz) = (xy)z\ \forall x,y \in \LL \}\\
\Nuc(\LL) &:=& \Nuc_l(\LL)\cap \Nuc_m(\LL)\cap \Nuc_r(\LL)
\end{array}
\]
Each of these is an associative subloop of $\LL$.\begin{lemma}
\label{lem:nuc-char}
If $\LL$ is a left loop, then
\[
\begin{array}{rcl}
L(\Nuc_l(\LL)) &=& \bigcap_{x\in \LL} L(\LL) L(x)\iv \\
L(\Nuc_m(\LL) &=& \bigcap_{x\in \LL} L(x)\iv L(\LL) .
\end{array}
\]
\end{lemma}

\begin{proof}
If $g\in \bigcap_{x\in\LL} L(\LL) L(x)\iv$, then $g = L(a)$ for some $a\in \LL$,
and for each $x\in\LL$, there exists $y\in\LL$ such that $L(y) = L(a)L(x)$. Applying both
sides to $1$ gives $y=a\cdot x$, and thus $L(a)L(x) = L(a\cdot x)$ for all $x\in\LL$, i.e.,
$a\in \Nuc_l(\LL)$. Reversing the argument yields the other inclusion, and the argument for
$\Nuc_m(\LL)$ is similar.
\end{proof}

Given a loop $\LL$, a subloop $\KK$ is said to be \emph{normal} if, for all $x,y\in \LL$,
$x\cdot (y\cdot \KK) = (x\cdot y)\cdot \KK$, $x\cdot \KK = \KK\cdot x$, and
$(\KK\cdot x)\cdot y = \KK\cdot (x\cdot y)$ (\cite{Br}, p. 60, IV.1). These three
conditions are clearly equivalent to the pair
\begin{equation}
\eqlabel{normality}
x\cdot (\KK\cdot y) = \KK\cdot (x\cdot y)
\qquad \text{and}\qquad
x\cdot (\KK\cdot y) = (x\cdot \KK)\cdot y
\end{equation}
for all $x,y\in \KK$.

\section{Twisted Subgroups}
\seclabel{ts}

Although the notion of a twisted subgroup of a group has been around for
some time (see Remark \remref{history}), we follow here
the terminology of Aschbacher \cite{Asch}, who proved one of the main
structural results about twisted subgroups (our Proposition \propref{Asch-main}
below). Our definition is a trivial modification of his.

For a subset $T$ of a group $G$, we use the notation
$T\iv := \{ x\iv : x\in T \}$ and $xTx := \{ xyx: y\in T \}$ for
$x\in T$.

\begin{definition}\cite{Asch}
\deflabel{ts}
A subset $T$ of a group $G$ is a \textit{twisted subgroup} of $G$ if
\item (i) $1 \in T$, \item(ii) $T\iv = T$, and \item(iii) $xTx\subseteq T$ for all $x \in T$.
\end{definition}

\begin{remark}
One can replace (ii) and (iii) with the equivalent assertion\hspace*{\fill}\linebreak
\noindent (ii$^{\prime}$) $xy\iv x \in T$ for all $x,y\in T$.
\end{remark}

A twisted subgroup $T$ of a group $G$ is said to be \emph{uniquely}
$2$-\emph{divisible} if each $x\in T$ has a unique square root in $T$, that is,
a unique element $x^{1/2}\in T$ such that $(x^{1/2})^2 = x$. As we do
with loops, we will abuse terminology slightly and drop the adverb ``uniquely".

An easy induction argument shows the following

\begin{proposition} (\cite{Asch}, Lem. 1.2(1))
\proplabel{mono}
Let $G$ be a group and let $T\subseteq G$ be a twisted subgroup.
Then for each $x\in T$, $\langle x \rangle \subseteq T$.
\end{proposition}

\begin{remark}
In some cases,  portions of Definition \defref{ts} are redundant.
\begin{enumerate}
\item[1.] For \textit{finite} groups, the proof of Proposition \ref{prop:mono}
shows that a subset satisfying (i) and (iii) necessarily satisfies (ii) (\cite{Asch}, Lemma 1.2(1)).
\item[2.] If $T\subseteq G$ is a left transversal of a subgroup $H\leq G$ such that
(iii) holds, then (ii) holds. Indeed, for $x\in T$, let $x^{\prime}\in T$
denote the representative of $x\iv H$. Then $x^{\prime} x x^{\prime}\in T$
and $x^{\prime} x x^{\prime}H = x^{\prime}H$, which implies
$x^{\prime} x x^{\prime} = x^{\prime}$. Thus $x^{\prime} = x\iv$.
In this case, the induced left loop $(T,\cdot)$ is a Bol loop; see
Proposition \propref{ts-bol}.
\item[3.] Glauberman showed that in the finite $2$-divisible
case, both (i) and (ii) are redundant (\cite{Gl1}, Lemma 3; \cite{Gl2}, Remark 7).
More precisely, he showed that if $T$ is a subset of a group $G$ satisfying (iii)
and such that every element of $T$ has finite odd order, then (i) and (ii) hold,
and every element of $T$ has a unique square root.
\end{enumerate}
\end{remark}

Of course, any subgroup is a twisted subgroup, but the notion of twisted subgroup
is modeled on the following example which is not a subgroup.

\begin{example}
\exlabel{1}
Let $G$ be a group, and fix $\tau\in \Aut(G)$. Define
\[
K(\tau) := \{ g\in G : g^{\tau} = g\iv \} .
\]
Then $K(\tau)$ is a twisted subgroup of $G$.
If $\tau^2 = 1$, define
\[
B(\tau) := \{ gg^{-\tau} : g\in G \}
\]
Then $B(\tau)$ is a twisted subgroup of $G$ and $B(\tau)\subseteq K(\tau)$.
\end{example}

\begin{example}
\exlabel{2}
Let $T$ be a twisted subgroup of a group $G$. For $x\in T$, define $\theta_x \in T!$ by
$\theta_x y = xyx$ for all $y\in T$. Then $\theta_1 = 1_{T!}$, $\theta_{x\iv} = \theta_x\iv$,
and $\theta_x \theta_y \theta_x = \theta_{xyx}$ for all $x,y\in T$. Thus $\hat{T} = \{ \theta_x : x\in T \}$
is a twisted subgroup of $T!$. For later reference, we will denote by $\hat{G}$ the subgroup of $T!$
generated by $\hat{T}$.
\end{example}

The \emph{associates} of a twisted subgroup $T$ of a group $G$
are the translates $aT = Ta\iv$, $a\in T$.

\begin{proposition} (\cite{Asch}, Lemma 1.5(1))
\proplabel{associates}
Every associate of a twisted subgroup is a twisted subgroup.
\end{proposition}

Most interesting results about twisted subgroups are predicated upon the assumption
that a twisted subgroup $T$ generates its group $G$. In this case we will just say that
$T$ is a generating twisted subgroup of $G$. Contained in such $T$ are important
normal subgroups of $G$. First we consider the intersection of all associates.

\begin{theorem}
\thmlabel{sharp}
Let $G$ be a group with generating twisted subgroup $T$, and let
\[
T^{\#} = \bigcap_{x\in T} xT .
\]
Then  $T^{\#} \subseteq T$,  $T^{\#} = \bigcap_{x\in T} Tx$, and $T^{\#} \triangleleft G$.
\end{theorem}

\begin{proof}
$T^{\#} \subseteq T$ is clear since $1\in T$, while $T^{\#} = \bigcap_{x\in T} Tx$
follows from $xT = Tx\iv$ for $x\in T$.
Now fix $a,b\in T^{\#}$ and $x\in T$. There exists $u,v\in T$ such that $a = xu \in xT$
and $b = u\iv v \in u\iv T$ (using $T\iv = T$). Hence $ab = xv \in xT$, and
since $x\in T$ is arbitrary, $ab \in T^{\#}$. Thus $T^{\#}$ is a subgroup of $G$. For
each $y\in T$,
\[
yT^{\#}y\iv = \bigcap_{x\in T} yxTy\iv= \bigcap_{x\in T} (yxy) T 
= \bigcap_{x\in T} \theta_y x T = T^{\#}.
\]
Since $T$ generates $G$, $T^{\#}$ is normal in $G$.
\end{proof}

A more important normal subgroup sitting inside a twisted subgroup was introduced by Aschbacher
(\cite{Asch}, p. 117). Our motivating discussion is a simplified version of his. Let $G$ be a group and
let $T\subseteq G$ be a generating twisted subgroup $G$. Consider the group
$G_0 = \langle ( x, x\iv ) : x\in T \rangle < G\times G$ generated by the graph
$\{ (x,x\iv) : x\in T \}$ of the inversion
mapping $x\mapsto x\iv$ on $T$. Let $\pi_i : G\times G\to G$ denote the projection onto the
$i$th factor. As a subgroup of $G\times G$, $G_0$ is invariant under the action of the swapping
automorphism $(x,y)\mapsto (y,x)$. This automorphism restricts to an isomorphism of the kernels
$\mathrm{Ker}(\pi_i |_{G_0})$. Each kernel is obviously isomorphic to the following subgroup of
$G$:
\[
T^{\prime} = \{ x_1 \cdots x_n : x_1\iv \cdots x_n\iv = 1,x_i \in T \}.\eqno{(1)}
\]
We have $T^{\prime} = \pi_1({\rm Ker}(\pi_2 |_{G_0})) \triangleleft \pi_1(G_0) = G$, since $T$
generates $G$. From the preceding discussion, we see that $G_0$ is the graph of an automorphism $\tau$
of $G$ if and only if  $T^{\prime} = \langle 1 \rangle$. In other words, $T$ is a subset of some $K(\tau)$
if and only if $T^{\prime} = \langle 1 \rangle$. This proves almost all of the following result.

\begin{proposition} (\cite{Asch}, Theorem 2.2)
\proplabel{Asch-main} 
Let $G$ be a group with generating twisted subgroup $T$. There
exists $\tau\in\Aut(G)$ with $\tau^2 = 1$ such that $T\subseteq K(\tau)$
if and only if $T^{\prime} = \langle 1 \rangle$. In this case the automorphism
$\tau$ uniquely determined.
\end{proposition}

\begin{proof}
All that remains is the uniqueness and the order. If $T\subseteq K(\sigma )$ for some
$\sigma\in\Aut(G)$, then $\tau\sigma\iv$ centralizes $T$, but then $\sigma = \tau$ since $T$
generates $G$. Since $T\subseteq K(\sigma )$ implies $T\subseteq K(\sigma\iv)$, it follows that $\tau^2 = 1$.
\end{proof}

For a twisted subgroup $T$ (whether it generates $G$ or not), we define its (Aschbacher)
\textit{radical} to be the normal subgroup $T^{\prime}$ given by (1) (\cite{Asch}, p.117).
If $T^{\prime} = \langle 1 \rangle$, then we say that $T$ is \textit{radical-free}.

If $T$ is radical-free and generates $G$, we will refer to the uniquely determined
$\tau\in\Aut(G)$ of order $2$ such that $T\subseteq K(\tau)$ as being the corresponding
\textit{Aschbacher automorphism}.

\begin{proposition}
Let $G$ be a group with generating twisted subgroup $T$.
Then $T^{\prime}$ is contained in every associate of $T$,
and hence $T^{\prime} \subseteq T^{\#}$.
\end{proposition}

\begin{proof}
The principal assertion is (\cite{Asch}, Theorem 2.1(3)), and the rest
follows from the definition of $T^{\#}$.
\end{proof}

\begin{remarks}\hspace*{\fill}
\begin{enumerate}
\item If $T$ is a \emph{proper} generating twisted subgroup of $G$, then
$T^{\prime}$ and $T^{\#}$ are proper normal subgroups of $G$. Thus twisted
subgroups of simple groups are radical-free and the intersection of all their
associates is trivial.
\item If $T$ is actually a subgroup of $G$, then the radical $T^{\prime}$
is the derived subgroup of $T$. This motivates our choice of notation, which is
different from that of \cite{Asch}.
\end{enumerate}
\end{remarks}
 
Besides the canonical projection of the previous proposition, there is another
radical-free twisted subgroup associated with any twisted subgroup. Here
we use the definitions and notation of Example \exref{2}.
 
\begin{theorem}
\thmlabel{hat-rad-free}
Let $G$ be a group, let $T\subseteq G$ be a twisted subgroup, and let
$\hat{T} = \{ \theta_x : x\in T \}$ and $\hat{G} = \langle \hat{T} \rangle$.
Then $\hat{T}$ is a radical-free twisted subgroup of $\hat{G}$.
\end{theorem}

\begin{proof}
If $\theta_{x_1}\cdots \theta_{x_n} \in \hat{T}^{\prime}$ for some $x_i \in T$, $1\le i\le n$,
then $1 = \theta_{x_n}\cdots \theta_{x_1} 1 = x_n\cdots x_1^2 \cdots x_n$, and rearranging
gives $x_1\cdots x_n^2 \cdots x_1 = 1$. But then for all $y\in T$,
\[
\theta_{x_1}\cdots \theta_{x_n} y
= \theta_{x_1}\cdots \theta_{x_n}\theta_{x_n}\cdots \theta_{x_1}y
=\theta_{x_1\cdots x_n^2 \cdots x_1} y = y .
\]
Thus $\theta_{x_1}\cdots \theta_{x_n} = 1_{\hat{G}}$, and therefore
$\hat{T}^{\prime} = \langle 1 \rangle$.
\end{proof}

In view of the preceding theorem, it is not surprising that the radical
is exactly the obstruction to a natural permutation representation of
a group $G$ on a generating twisted subgroup $T$.

\begin{theorem}
\thmlabel{extend}
Let $G$ be a group with generating twisted subgroup $T$. The mapping
$\theta : T\to \hat{T}$ defined by $\theta_x y = xyx$ ($x,y\in T$)
extends to a homomorphism $\theta : G\to \hat{G}$ if and only if
$T^{\prime} = \langle 1\rangle$. In this case,
$\Ker (\theta) = Z(G)\cap C_G(\tau)$ where $\tau\in\Aut(G)$ is the
Aschbacher automorphism.
\end{theorem}

\begin{proof}
The subgroup $\langle (x,\theta_{x}) : x\in T \rangle$ of $G\times T!$
is the graph of a homomorphism from $G$ into $T!$ if and only if the group
$K = \{ \theta_{x_1} \cdots \theta_{x_n} : x_1\cdots x_n = 1, x_i\in T \}
= \langle 1 \rangle$. The mapping $K\to T; \phi\mapsto \phi 1$ is
a homomorphism with image in $T^{\prime}$.
Indeed, if $\phi = \theta_{x_1}\cdots \theta_{x_n}$ for $x_i\in T$ with
$x_1\cdots x_n = 1$, then $\phi 1 = x_n\cdots x_1 \in T^{\prime}$.
This homomorphism is clearly onto, and if $\phi 1 = 1$, then
$x_n\cdots x_1 = 1$, whence
$\phi = \theta_{x_1}\cdots \theta_{x_n} = 1_{\hat{G}}$.
Therefore $K$ is isomorphic to $T^{\prime}$. This establishes the
first assertion. Assume now that $T^{\prime} = \langle 1\rangle$ and let
$\tau\in\Aut(G)$ denote the Aschbacher automorphism.
Fix $g = x_1\cdots x_n \in \Ker (\theta)$ for $x_i\in T$. Then
$g y x_n\cdots x_1 = y$ for all $y\in T$. Taking $y=1$, we have
$x_n\cdots x_1 = g\iv$. Thus $g$ centralizes $T$. Since $T$
generates $G$, $g\in Z(G)$. Also,
$g^{\tau} = x_1^{\tau}\cdots x_n^{\tau} = (x_n\cdots x_1)\iv = g$.
Conversely, if $g = x_1\cdots x_n \in Z(G)\cap C_G(\tau)$ for
$x_i\in T$, then
$g\iv = g^{-\tau} = x_n^{-\tau}\cdots x_1^{-\tau} = x_n\cdots x_1$,
and so $g\in \Ker (\theta)$.
\end{proof}

\begin{corollary}
\corolabel{ts-extend}
Let $G$ be a simple group with generating twisted subgroup $T$. The
mapping $\theta : T\to \hat{T}$ defined by $\theta_x y = xyx$ ($x,y\in T$)
extends to an isomorphism $\theta : G\to \hat{G}$.
\end{corollary}

\section{$2$-Divisible Twisted Subgroups and B-loops}
\seclabel{2dts}

We now focus our attention on $2$-divisible twisted subgroups and their 
associated B-loops. Much (though not all) of this section is an adumbration of
Glauberman's fundamental results \cite{Gl1, Gl2}. We give (often
simpler) proofs of some of his results to make the exposition self-contained.

\begin{lemma}
\lemlabel{equiv}
Let $G$ be a group, and let $T\subseteq G$ be a twisted subgroup in which every element
has finite order. The following are equivalent: (i) $T$ is $2$-divisible; (ii) every element of $T$
has odd order; (iii) no element of $T$ has order $2$. If, in addition, $T$ has finite order, then
these conditions imply: (iv) $|T|$ is odd.
\end{lemma}

\begin{proof}
If $T$ is $2$-divisible, then obviously no elements of $T$ have even order, and so (i) implies
(ii) and (iii). The equivalence of (ii) and (iii) is trivial. Now assume (ii) and let $z\in T$ be given
with order $2k+1$. Then $z^{k+1}$ is a square root of $z$ in $T$. If $y\in T$ were another
square root of $z$, then $y^{2(2k+1)} = 1$, and so $y^{2k+1} = 1$. But then
$z^{2k+1} = y^{2k+1} = z^k y$ so that $y = z^{k+1}$. Thus (i) holds. For the remaining
assertion, note that the inversion mapping $x\mapsto x\iv$ is a permutation of the set
$T \backslash \{ 1 \}$. If $T$ had even order, then this mapping would necessarily fix some
$a \neq 1$. But then $a^2 = 1$, whence $T$ is not $2$-divisible. Thus (i) implies (iv).
\end{proof}

\begin{example}
\exlabel{converse}
In general, condition (iv) of Lemma \lemref{equiv} does not imply the other
conditions. Indeed, let $G = S_3$, the symmetric group on $3$ letters and
let $T$ be the set of transpositions. Then $|T| = 3$, but every element of $T$
has order $2$.
\end{example}

Radical-free, generating, $2$-divisible twisted subgroups are ``rigid"
in the sense that they are uniquely determined by the Aschbacher
automorphism. For a subset $S$ of a group $G$, we denote
$S^2 = \{ x^2 : x\in S\}$.

\begin{theorem}
\thmlabel{rigid}
Let $G$ be a group, let $T\subseteq G$ be a radical-free,
generating twisted subgroup, and let $\tau\in\Aut(G)$
be the corresponding Aschbacher automorphism. Then
\[
K(\tau)^2 \subseteq B(\tau) \subseteq T \subseteq K(\tau).
\]
In particular, if $T$ is $2$-divisible, then $B(\tau) = T = K(\tau)$,
and $T$ is a left transversal of $C_G(\tau)$ in $G$.
\end{theorem}\begin{proof} For $g\in K(\tau)$, $g^2 = gg^{-\tau} \in B(\tau)$, and
so $K(\tau)^2 \subseteq B(\tau)$. For $g\in G$, $g = x_1\cdots x_n$ for
some $x_i\in T$. Thus $gg^{-\tau} = x_1\cdots x_n x_n\cdots x_1 \in T$,
since $T$ is a twisted subgroup. Thus $B(\tau) \subseteq T$. The equality
in the $2$-divisible case follows immediately. Now for $g\in G$, we have
$a := (gg^{-\tau})^{1/2}\in T$, and it is easy to check that
$h := a\iv g \in C_G(\tau)$. The uniqueness of the decomposition
$g = ah$ is obvious.
\end{proof}

\begin{definition}
\label{defn:B-loop}
Let $G$ be a group with a $2$-divisible twisted subgroup $T$.
Define a binary operation $\odot : T\times T\to T$ by 
\[
x\odot y := (x y^2 x)^{1/2}
\]
for $x,y\in T$. We follow Glauberman's notation \cite{Gl2} and denote the
magma $(T,\odot)$ by $T(1/2)$. We denote
the left multiplication maps for $T(1/2)$ by $b_x y := (xy^2 x)^{1/2}$ for $x,y\in T$.
\end{definition}

\begin{lemma} 
\label{lem:B-loop}
Let $G$ be a group with a $2$-divisible twisted subgroup $T$.
\begin{enumerate}
\item[1.] $T(1/2)$ is a B-loop.
\item[2.] Integer powers of elements in $T$ formed in $G$ agree with those in $T(1/2)$.
Thus an element has finite order in $T$ if and only if it has the same order in $T(1/2)$.
\item[3.] If $T$ is radical-free and generates $G$, then $T(1/2)$ agrees with the left loop
structure induced on $T$ as a left transversal.
\end{enumerate}
\end{lemma}

\begin{proof}
(1) and (2) follow from Lemma 3 in \cite{Gl1} and the following remark. For (3), 
let $\tau\in\Aut(G)$ be the Aschbacher automorphism, and note that for $x, y\in T$,
$x\odot y = ((xy)(xy)^{-\tau})^{1/2}$.
\end{proof}

\begin{remark}
It is slightly more common in the loop theory literature to use the operation
\[
x\odot^{\prime} y := x^{1/2}yx^{1/2}
\]
for $x,y\in T$. Clearly the squaring map $x\mapsto x^2$ is an isomorphism of $(T,\odot )$ onto
$(T,\odot^{\prime})$. That some authors prefer $\odot^{\prime}$ is partly because the B-loop
$(T,\odot^{\prime})$ is isotopic to a quasigroup structure on $T$ given by $(x,y)\mapsto xy\iv x$.
(For the notion of isotopy, see any of the standard references \cite{Bel2, Br, Pf}.)
With different terminology than that used here, the preceding construction on $2$-divisible twisted subgroups
(using either $\odot$ or $\odot^{\prime}$) can be found in Foguel and Ungar \cite{FU},
Glauberman (\cite{Gl1}, Lemma 3), Kikkawa (\cite{Kik}, Theorem 5), Kreuzer \cite{Kr}, and
Kiechle (\cite{Kie}, Chap. 6D). There is
a related construction in uniquely $2$-divisible loops $\LL$ which goes as follows: for $x,y\in\LL$, define
$x\odot^{\prime\prime} y = x^{1/2} \cdot (y\cdot x^{1/2})$. For certain loops $(\LL,\cdot )$, the new magma
$(\LL, \odot^{\prime\prime})$ turns out to be a B-loop. For $2$-divisible Moufang loops, this construction is due
to Bruck (\cite{Br}, VII.5.2, p. 121). For uniquely $2$-divisible Bol loops, it is implicit in the work of Belousov
(\cite{Bel1}, \cite{Bel2}), and is spelled out in work of P.T. Nagy and K. Strambach (\cite{NS}, p. 301, Thm. 7)
as well as the recent dissertation of G. Nagy (\cite{Na3}, T\'{e}tel 2.3.6). As it turns out, these loop-based
constructions are no more general than the construction for twisted subgroups, because they all depend on
the fact that the loop in question can be identified in a natural way with a twisted subgroup. The B-loop structure
is then transferred from the twisted subgroup to the loop. We will see how this works for Bol loops in
\S\secref{2dbol}.
\end{remark}

The following is clear from the definitions.

\begin{lemma}
\lemlabel{conj}
Let $G$ be a group with $2$-divisible twisted subgroup $T$,
and let $s \in T!$ denote the squaring map on T: $s(x) = x^2$.
Then $b_x = s\iv\theta_x s$ for all $x\in T$. Thus 
$\LMlt(T(1/2))$ is conjugate in $T!$ to
$\hat{G} = \langle \theta_x : x\in T\rangle$.
\end{lemma}

\begin{theorem}
\thmlabel{b-extend}
Let $G$ be a group with a generating, $2$-divisible twisted 
subgroup $T$. The mapping $b: T\to T!$ defined by $b_x y = (xy^2x)^{1/2}$
($x,y\in T$) extends to a homomorphism $b : G\to T!$ if and only if
$T^{\prime} = \langle 1\rangle$. In this case, there is an exact sequence
\[
1\to Z(G)\cap C_G(\tau) \to G \to \LMlt(T(1/2)) \to 1
\]
where $\tau\in\Aut(G)$ is the Aschbacher automorphism, and
$Z(G)\cap C_G(\tau)$ is the core of $C_G(\tau)$ in $G$.
\end{theorem}

\begin{proof}
This follows from Lemma \lemref{conj}, Theorem \thmref{extend}
and Proposition \propref{core}.
\end{proof}

One reason it is particularly convenient to work with the B-loop associated
with a $2$-divisible twisted subgroup is the following.

\begin{lemma} (cf. \cite{Gl1}, p. 379, Lemma 4)
\lemlabel{subloops}
Let $G$ be a group and let $T\subseteq G$ be a $2$-divisible twisted
subgroup. Then $K\subseteq T$ is a twisted subgroup of $G$ if and
only if $K(1/2) := (K,\odot)$ is a subloop of $T(1/2)$.
\end{lemma}

\begin{proof}
This is immediate from the definition of $\odot$.
\end{proof}

The second corollary to the following result is (\cite{Gl1}, p. 384, Corollary 3).

\begin{theorem}
\thmlabel{preLagrange}
Let $G$ be a finite group, let $T\subseteq G$ be a $2$-divisible
twisted subgroup, and let $A\subseteq T$ be a subgroup of $G$.
\begin{enumerate}
\item[1.] If $A$ is normal in $G$, then $|A|$ divides $|T|$.
\item[2.] If $A$ is abelian, then $|A|$ divides $|T|$.
\end{enumerate}
\end{theorem}

\begin{proof}
(1) For each $x\in T$, note that $x A = x^{1/2} A x^{1/2} \subseteq T$,
and thus $\{ xA : x\in T \}$ partitions $T$ into subsets of equal cardinality.

(2) $A(1/2)$ is an abelian group isomorphic to $A$. 
The restriction of 
\[
b: T\to \LMlt(T(1/2)); x \mapsto (y \mapsto x\odot y)
\]
to $A$ is a homomorphism of $A(1/2)$ onto its image. The orbits
$\{ \{ b_x y : x\in A \} : y\in T\}$ clearly partition $T$, and the orbit
through $1\in T$ is $A$ itself since $A$ is $2$-divisible. The action of
$A$ on any orbit is regular since $T(1/2)$ is a loop.
\end{proof}

\begin{corollary}
\corolabel{rad-divides}
Let $G$ be a group and let $T\subseteq G$ be a finite $2$-divisible
twisted subgroup. Then $|T^{\#}|$ and $|T^{\prime}|$ divide $|T|$.
\end{corollary}

\begin{corollary}
(Lagrange's Theorem)
\corolabel{lagrange}
Let $G$ be a group and let $T\subseteq G$ be a finite $2$-divisible
twisted subgroup. Then for every $x\in T$, the order of $x$ divides $|T|$.
\end{corollary}

\begin{proof}
Apply Theorem \ref{thm:preLagrange}(1) to the subgroup
$\langle x\rangle \subseteq T$.
\end{proof}

\begin{remark}
\remlabel{lagrange-no}
As Example \exref{converse} indicates, Lagrange's theorem does
not hold for all twisted subgroups.
\end{remark}

The following is a distilled version of (\cite{Gl2}, Theorem 14).

\begin{theorem}
\thmlabel{odd}
Let $G$ be a finite group, and let $T\subseteq G$ be a $2$-divisible, generating twisted subgroup.
Then $G$ has odd order.
\end{theorem}

\begin{proof}
Assume first that $T$ is radical-free and let $\tau\in\Aut(G)$ denote the
Aschbacher automorphism. By Theorem \thmref{rigid}, $T = B(\tau)$.
By Glauberman's $Z^*$ Theorem (\cite{Gl3}, Theorem 1), there exists a normal
subgroup $N$ of $G\langle \tau \rangle$ such that $|N|$ is odd and 
$\tau N \in Z(G\langle \tau \rangle / N)$. But then for all $g\in G$, 
$g \tau g\iv \tau = gg^{-\tau} \in N$. Thus $T\subseteq N$. Since $T$
generates $G$, $G = N$. For the general case, $G/T^{\prime}$ must have
odd order, and thus by Corollary \cororef{rad-divides},
$|G| = |G/T^{\prime}| |T^{\prime}|$ is odd.
\end{proof}

\begin{definition}
\deflabel{pi}
Let $\pi$ be a set of primes. A positive integer $n$ is a $\pi$-\emph{number} if $n=1$
or if $n$ is a product of primes in $\pi$. For every positive integer $n$, let $n_\pi$ denote
the largest $\pi$-number that divides $n$.  As usual, a finite group $G$ is a
$\pi$-\emph{group} if $|G| = |G|_{\pi}$. If $T\subseteq G$ is a twisted subgroup, then
we say that $T$ is a \emph{twisted} $\pi$-\emph{subgroup} of $G$ if $|T| = |T|_{\pi}$.
We say that $T$ satisfies the \emph{Hall} $\pi$-\emph{condition} if there exists a twisted
$\pi$-subgroup $S$ of $G$ such that $S\subset T$ and $|S| = |T|_{\pi}$. If
$\pi = \{ p\}$, we say that $T$ satisfies the \emph{Sylow} $p$-\emph{condition} if
$T$ satisfies the Hall $\{p\}$-condition.
\end{definition}

\begin{lemma}
\lemlabel{preHall}
Let $G$ be a finite group of odd order, let $\pi$ be a set of primes, and let
$\beta\in\Aut(G)$ have order $2$. Then every $\pi$-subgroup of $G$ fixed
by $\beta$ is contained in a Hall $\pi$-subgroup of $G$ fixed by $\beta$.
\end{lemma}

\begin{proof} 
Since $G$ is solvable \cite{FT}, this is just (\cite{Gl1}, p. 391, Lemma 11).
\end{proof}

Glauberman remarked that the following result can be established by
purely group-theoretical means (\cite{Gl2}, p. 413, Remark 7).

\begin{theorem} (\cite{Gl2}, Theorem 15)
\thmlabel{pi}
Let $G$ be a finite group, let $T\subseteq G$
be a $2$-divisible, generating twisted subgroup, and 
let $\pi$ be a set of odd primes. Then $T$ is a twisted
$\pi$-subgroup if and only if $G$ is a $\pi$-group.
\end{theorem}

\begin{proof}
By Theorem \thmref{rigid}, $|T|$ divides $|G|$, and so if
$G$ is a $\pi$-group, then $T$ is certainly a twisted 
$\pi$-subgroup. For the converse, assume first that $T$
is radical-free, let $\tau\in\Aut(G)$ denote the Aschbacher automorphism,
and set $H := C_G(\tau)$. Let $P_0$ be a Hall $\pi$-subgroup
of $H$. By Theorem \thmref{odd}, $|G|$ is odd, and so by
Lemma \lemref{preHall}, $P_0$ is contained in some Hall $\pi$-subgroup
$P$ of $G$ which is fixed by $\tau$.  By Theorem \thmref{rigid},
$S := P\cap T$ is a left transversal of $P_0 = P\cap H$ in $P$, and hence
$| S | = |P| / |P_0| = |G|_{\pi} / |H|_{\pi} = [ G : H ]_{\pi} = |T|_{\pi} = |T|$.
Thus $T\subseteq P$, and since $T$ generates $G$, we have $G = P$.
In the general case, $G/T^{\prime}$ is a $\pi$-group, and so by Corollary
\cororef{rad-divides}, $|G| = |G/T^{\prime}| |T^{\prime}|$ is a $\pi$-number.
\end{proof}

\begin{theorem}
(Hall's Theorem, cf. \cite{Gl1}, p. 392, Theorem 8)
\thmlabel{hall}
Let $G$ be a finite group, and let $T\subseteq G$ be a $2$-divisible twisted
subgroup. For every set $\pi$ of primes, $T$ satisfies
the Hall $\pi$-condition, and thus $T(1/2)$ has a Hall $\pi$-subloop.
\end{theorem}

\begin{proof}
Without loss of generality, assume $T$ generates $G$. Since $T$ can be identified
with its radical-free image $b(T)\subseteq \LMlt(T(1/2))$, there is no loss of generality
in assuming that $T$ is radical-free. Repeating the proof of Theorem \thmref{pi},
we obtain a Hall $\pi$-subgroup $P$ of $G$ fixed by $\tau$ such that
$S := P\cap T$ satisfies $|S| = |T|_{\pi}$. $S(1/2)$ is a Hall $\pi$-subloop
of $T(1/2)$ by Lemma \lemref{subloops}.
\end{proof}

\begin{remark}
\remlabel{hall}
Using Theorem \thmref{pi} (i.e., \cite{Gl2}, Theorem 15), 
Glauberman showed that the Hall $\pi$-subloops of $T(1/2)$
are all conjugate under $C_G(\tau)$, that every prime dividing the number
of such subloops also divides $|T|$ and is not in $\pi$, and that every
$\pi$-subloop of $T(1/2)$ (that is, every twisted $\pi$-subgroup of
$G$ contained in $T$) is contained in a Hall $\pi$-subloop; see (\cite{Gl1}, Theorem 8).
\end{remark}

\begin{corollary}
\corolabel{sylow}
(Sylow's Theorem, \cite{Gl1}, p. 394, Corollary 3)
Let $G$ be a finite group with a $2$-divisible twisted subgroup $T$.
For every prime $p$, $T$ satisfies the Sylow $p$-condition, and thus
$T(1/2)$ has a Sylow $p$-subloop.
\end{corollary}

\begin{remark}
In \cite{Gl1}, Glauberman originally gave separate proofs under different
hypotheses of the Sylow and Hall theorems for B-loops,
because at the time it was not known that the group generated by a
$2$-divisible twisted subgroup must have odd order. In light of his later result
(\cite{Gl2}, Theorem 14), our Theorem \thmref{odd}, the Sylow
result easily follows from the Hall result. The additional properties
mentioned in Remark \remref{hall} obviously hold in the Sylow case
as well.
\end{remark}

\begin{corollary}
(Cauchy's Theorem, cf. \cite{Gl1}, p. 394, Corollary 1)
\corolabel{ts-cauchy}
Let $G$ be a finite group, and let $T\subseteq G$ be a $2$-divisible twisted subgroup.
If a prime $p\mid |T|$, then $T$ contains an element of order $p$.
\end{corollary}

\begin{proof}
By Corollary \cororef{sylow}, there is a twisted $p$-subgroup $S$ of $G$
such that $S\subseteq T$ and $|S| = |T|_p$. The rest follows from Lagrange's theorem
(Corollary \cororef{lagrange}).
\end{proof}

\begin{remark}
\remlabel{no-hope}
There is no hope of extending Cauchy's theorem to all twisted subgroups.
There exists a twisted subgroup of order $180$ which does not have an
element of order $5$. We will discuss this further in
Remark \remref{bol-no-hope}
\end{remark}

\begin{proposition}
(Strong Lagrange Theorem)
\proplabel{strongLagrange}
Let $G$ be a finite group, and let $T\subseteq G$ be a $2$-divisible
twisted subgroup. If $A \subseteq B \subseteq T$ are
twisted subgroups of $G$, then $|A|$ divides $|B|$.
\end{proposition}

\begin{proof}
(\cite{Gl1}, p. 395, Corollary 4).
\end{proof}

\begin{remark}
Feder \cite{Feder} recently extended Proposition \propref{strongLagrange}
to \emph{strong near subgroups}, which include twisted subgroups of
odd order as a special case. Roughly speaking, strong near subgroups
are twisted subgroups in which the $2$-elements are well-behaved.
\end{remark}

\section{Bol loops}
\seclabel{bol}

We now apply the results of \S\secref{ts} to Bol loops. In fact, Bol loops are
related to twisted subgroups in more than one way.

\begin{example}
\exlabel{3}
Let $\LL$ be a loop, and let $L(\LL) = \{ L(x) : x\in \LL \}$ denote its set of left translations.
Then $\LL$ is a Bol loop if and only $L(\LL)$ is a twisted subgroup of $\LMlt(\LL)$. Also,
if $\KK$ is a subloop of $\LL$, then $L(\KK)$ is a twisted subgroup of $\LMlt(\LL)$.
\end{example}

More generally, we have the following.

\begin{proposition} (\cite{KJ}, Remark 4.4(2))
\proplabel{ts-bol}
Let $G$ be a group, $H\leq G$, and $T\subseteq G$ a transversal of $H$.
If $T$ is a twisted subgroup, then $(T,\cdot)$ is a Bol loop. Conversely,
if $H$ is core-free and $(T,\cdot)$ is a Bol loop, then $T$ is a twisted subgroup.
\end{proposition}

\begin{example}
\exlabel{4}
Let $\LL$ be a Bol loop. For each $x\in\LL$, set $P(x) = L(x)R(x)$, and let
$P(\LL) = \{ P(x) : x\in \LL \}$. Then $P(\LL)$ is a twisted subgroup of the
group $\PMlt(\LL) := \langle P(x) : x\in\LL \rangle$. This is really just a special
case of Example \exref{2}. Indeed, for $x,y\in \LL$, we have
\begin{equation}
\eqlabel{PL}
\theta_{L(x)} L(y) = L(P(x)y).
\end{equation}
Thus for $x,y,z\in \LL$, we compute
\[
\begin{array}{rcl}
L(P(x\cdot (y\cdot x))z) &=& \theta_{L(x\cdot (y\cdot x))} L(z) = \theta_{L(x)L(y)L(x)}L(z) \\
&=& \theta_{L(x)} \theta_{L(y)} \theta_{L(x)} L(z) = L(P(x)P(y)P(x)z).
\end{array}
\]
Thus $P(x\cdot (y\cdot x)) = P(x)P(y)P(x)$ as claimed. The other properties of twisted
subgroups follow similarly.
\end{example}

\begin{example}
\exlabel{5}
Let $\LL$ be a Bol loop. Then for each $x\in\LL$, the triple
\[
B(x) = (P(x),L(x\iv),L(x))
\]
is an autotopism of $\LL$. (For the notion of autotopism, see any of the standard
references \cite{Bel2, Br, Pf}.) Conversely, if $\LL$ is a loop in which each $B(x)$ is an autotopism,
then $\LL$ is a Bol loop. Let $\Btp(\LL) = \langle B(x): x\in\LL\rangle$ denote the group of all
\textit{Bol autotopisms} of $\LL$. Then from Examples~\exref{3} and \exref{4}, we see that the
set $B(\LL) = \{ B(x): x\in\LL \}$ is a twisted subgroup of $\Btp(\LL)$ (or of the entire autotopism
group of $\LL$). Geometrically, the Bol autotopism group $\Btp(\LL)$ is isomorphic to a subgroup
of the collineation group of the associated $3$-net, namely the direction-preserving collineation group
generated by Bol reflections \cite{FN}.
\end{example}

Recall that for a group $G$ with twisted subgroup $T$, the group $\hat{G}\subseteq T!$
is defined by $\hat{G} = \langle \theta_x : x\in T\rangle$; see Example \exref{2}.

\begin{lemma}
\lemlabel{Piso}
Let $\LL$ be a Bol loop. Then $\PMlt(\LL) \cong \widehat{\LMlt}(\LL)$.
The isomorphism is defined on generators by $P(x)\mapsto \theta(L(x))$.
In case $\LL$ is $2$-divisible, we also have $\PMlt(\LL) \cong \LMlt(\LL(1/2))$
\end{lemma}

\begin{proof}
The first assertion follows from \peqref{PL} in Example \exref{4}. The
second follows from Lemma \lemref{conj}.
\end{proof}

The distinction, therefore, between $\PMlt(\LL)$ and $\widehat{\LMlt}(\LL)$
is that the former acts directly on the loop $\LL$, while the latter acts on
the transversal $L(\LL)$.

\begin{corollary}
\corolabel{P-radfree}
Let $\LL$ be a Bol loop. Then $P(\LL)$ is a radical-free twisted
subgroup of $\PMlt(\LL)$.
\end{corollary}

\begin{proof}
This follows from Lemma \lemref{Piso} and Theorem \thmref{hat-rad-free}.
\end{proof}

First we consider the interpretation in $\LL$ of the normal subgroup
$T^{\#}$ for $T = \LMlt(\LL)$.

\begin{theorem}
\thmlabel{leftnuc}
If $\LL$ is a Bol loop, then $L(\LL)^{\#} = L(\Nuc_l(\LL)) = L(\Nuc_m(\LL))$.
\end{theorem}

\begin{proof}
$\LL$ has LIP, and so by Lemma \lemref{nuc-char},
$L(\Nuc_m(\LL)) = \bigcap_{x\in\LL} L(x\iv)L(\LL) = L(\LL)^{\#}$.
The other equality follows Theorem \thmref{sharp}.
\end{proof}

\begin{remark}
The equality $\Nuc_l(\LL) = \Nuc_m(\LL)$ for left loops with LIP is well-known
(e.g., \cite{Kie}, p. 62, (5.7)). Expressed in terms of a subset $T$ (such as
$L(\LL)$) of a group (such as $\LMlt(\LL)$), this just says that the equality
$\bigcap_{x\in T} xT = \bigcap_{x\in T} Tx$ holds provided that $T\iv = T$.
\end{remark}

\begin{corollary} (\cite{Na2}, p. 405, Lemma 1)
Let $\LL$ be a Bol loop. Then $\Nuc_l(\LL)$ is a normal subloop.
\end{corollary}

\begin{proof}
Using Theorems \ref{thm:leftnuc} and \ref{thm:sharp}, the conditions
of \peqref{normality} are easily checked.
\end{proof}

Next we turn to the radical.

\begin{definition}
\deflabel{bol-rad}
Let $\LL$ be a Bol loop. The \emph{radical} of $\LL$ is the set
$\LL^{\prime} := \{ x\in \LL  : L(x) \in L(\LL)^{\prime} \}$.
In case $\LL^{\prime} = \{ 1\}$ (i.e., $L(\LL)^{\prime} = \langle 1\rangle$),
we will say that $\LL$ is \textit{radical-free}.
\end{definition}

\begin{corollary}
\corolabel{rad-assoc}
Let $\LL$ be a Bol loop with radical $\LL^{\prime}$. Then $\LL^{\prime}$ is
an associative normal subloop of $\LL$ contained in $\Nuc_l(\LL)$.
\end{corollary}

\begin{proof}
That $\LL^{\prime} \subseteq \Nuc_l(\LL)$ follows from
$L(\LL)^{\prime}\subseteq L(\LL)^{\#}$ and Theorem \ref{thm:leftnuc}.
Thus $\LL^{\prime}$ is associative and the conditions of
\peqref{normality} follow easily. 
\end{proof}

\begin{remarks}\hspace*{\fill}
\label{rem:bol-rads}
\begin{enumerate}
\item As isomorphic abstract groups, there is, of course, no meaningful distinction between
the radical $\LL^{\prime}$ of a Bol loop $\LL$ and the radical $L(\LL)^{\prime}$
of the twisted subgroup $L(\LL)$ of left translations, particularly if one identifies the
loop $\LL$ with the induced loop structure on the transversal $L(\LL)$. However,
the distinction does help clarify the normality of $\LL^{\prime}$ as a sub\textit{loop}
of $\LL$ versus the normality of $L(\LL)^{\prime}$ as a \textit{subgroup} of $\LMlt(\LL)$.
\item Let $\LL$ be a Bruck loop, and let $\tau\in\Aut(\LMlt(\LL))$ denote conjugation by
the inversion mapping $x\mapsto x\iv$. Then the automorphic inverse property is equivalent
to $L(x)^{\tau} = L(x\iv)$ for all $x\in \LL$. Thus $\tau$ is the Aschbacher automorphism
of $\LMlt(\LL)$, and hence $\LL$ is a radical-free Bol loop.
\end{enumerate}
\end{remarks}

\begin{theorem}
\thmlabel{rep1}
Let $\LL$ be a Bol loop, and let $G = \LMlt(\LL)$.
The mapping $L(\LL)\to P(\LL); L(x)\mapsto P(x)$
extends to a homomorphism from $G$ onto $\PMlt(\LL)$ if and
only if $\LL^{\prime} = \{ 1\}$. In this case, the kernel of the
homomorphism is $Z(G)\cap C_G(\tau)$ where $\tau\in\Aut(G)$
is the Aschbacher automorphism. 
\end{theorem}

\begin{proof}
This follows from Lemma \lemref{Piso} and Theorem \thmref{extend}.
\end{proof}

Next we consider the Bol autotopism group $\Btp(\LL)$ of a Bol
loop $\LL$. Let $\Phi_i : \Btp(\LL)\to \Mlt(\LL); (f_1, f_2, f_3)\mapsto f_i$
denote the projection onto the $i^{\text{th}}$ component. Clearly
$\Phi_1$ is an epimorphism onto $\PMlt(\LL)$ and $\Phi_2$ and $\Phi_3$
are epimorphisms onto $\LMlt(\LL)$.

In the Bol loop context, the subloop we call the radical made
its first appearance in work of M. Funk and P. Nagy (\cite{FN},
p. 67, Theorem 1). The following is the algebraic version of their
geometric result.

\begin{theorem}
\thmlabel{btp-rad}
Let $\LL$ be a Bol loop, and let $\Phi_3 : \Btp(\LL) \to \LMlt(\LL)$ be the
projection onto the third factor. Then $\Ker(\Phi_3) \cong \LL^{\prime}$.
\end{theorem}

\begin{proof}
$(f,g,1)\in \Ker(\Phi_3)$ if and only if $g$ can be written as
$g = L(x_1\iv)\cdots L(x_n\iv)$ for some $x_i\in\LL$,
$i=1,\ldots ,n$, such that $L(x_1)\cdots L(x_n) = I$. Thus
$(f,g,1)\in \ker(\Phi_3)$ if and only if $g \in L(\LL)^{\prime}$,
and so the restriction of $\Phi_2$ to $\Ker(\Phi_3)$ is an isomorphism
onto $L(\LL)^{\prime}$.
\end{proof}

\begin{remark}
In particular, if $\LL$ is a radical-free Bol loop, then the group
$\Btp(\LL)$ simultaneously encodes both the graph of the 
Aschbacher automorphism $\tau\in\Aut(\LMlt(\LL))$ and the
graph of the homomorphism $\LMlt(\LL)\to \PMlt(\LL)$ described
in Theorem \thmref{rep1}. These are given by, respectively,
$f_3\mapsto f_2$ and $f_3\mapsto f_1$ for $(f_1,f_2,f_3)\in \Btp(\LL)$.
\end{remark}

\begin{lemma}
\lemlabel{kernel}
Let $\LL$ be a loop and set $G := \LMlt(\LL)$. 
Then $Z(G) = G\cap \{ R(x) : x \in \Nuc_r(\LL)\}$.
Therefore the set $\MM := \{ x\in \Nuc_r(\LL) : R(x)\in G\}$
is an abelian group.
\end{lemma}

\begin{proof}
An element $a\in\LL$ is in $\Nuc_r(\LL)$ if and only
if $R(a)$ centralizes $G$ in the full multiplication group
$\Mlt(\LL)$. So if some such $R(a)\in G$, then
$R(a)\in Z(G)$. Conversely, if $g \in Z(G)$, then setting
$a = g1$, we have $x\cdot a = L(x)g1 = gL(x)1 = gx$,
and so $g = R(a)$ and $a\in \Nuc_r(\LL)$. The rest
follows because the mapping $R : \MM\to Z(G);
x\mapsto R(x)$ is an anti-isomorphism.
\end{proof}

\begin{theorem}
Let $\LL$ be a Bol loop, let $G = \LMlt(\LL)$, and let
$\Phi_1 : \Btp(\LL) \to \PMlt(\LL)$ be the projection onto
the first factor. Then $\Ker(\Phi_1) \cong Z(G)\cap \{ g :
g = L(x_1)\cdots L(x_n) = L(x_1\iv)\cdots L(x_n\iv),
x_i\in \LL, i=1,\ldots,n\}$. If $\LL$ is radical-free, then
$\Ker(\Phi_1) \cong Z(G)\cap C_G(\tau)$ where $\tau\in\Aut(G)$
is the Aschbacher automorphism. 
\end{theorem}

\begin{proof}
A triple $(1,f_2,f_3)$ of permutations is an autotopism if and only
if $f_2 = f_3 = R(a)$ where $a = f_2(1) \in \Nuc_r(\LL)$.
As in the proof of Lemma \ref{lem:kernel}, this holds if and only if
$R(a)$ centralizes $\LMlt(\LL)$ in $\Mlt(\LL)$, and so 
$(1,R(a),R(a))\in \Btp(\LL)$ if and only if $R(a)\in Z(\LMlt(\LL))$
and $R(a) = L(x_1)\cdots L(x_n) = L(x_1\iv)\cdots L(x_n\iv)$ for
some $x_i\in\LL, i=1,\ldots,n\}$. The remaining assertion follows
immediately. 
\end{proof}

\begin{corollary}
Let $\LL$ be a radical-free Bol loop, let $G = \LMlt(\LL)$, and let
$\tau\in\Aut(G)$ be the Aschbacher automorphism. If
$Z(G)\cap C_G(\tau)  = \langle 1\rangle$,  then
$\Btp(\LL) \cong \LMlt(\LL) \cong \PMlt(\LL)$.
\end{corollary}

\begin{remark}
\remlabel{history}
Before proceeding on to $2$-divisible Bol loops, it is probably worthwhile to
insert a few historical remarks. The concept of twisted subgroup (though obviously
not the terminology we have adopted), and its relationship with quasigroup and loop theory, 
has been around for some time, and is not limited to the connection with Bol loops. For
example, a \textit{Fischer group} is a group $G$ and a subset $T\subseteq G$ of involutions
which generate $G$ such that for all $x,y\in T$, $(xy)^3 = 1$, and $xyx \in T$. If
$1\in T$, then $T$ is a twisted subgroup. Fischer groups arise in the study of distributive,
symmetric quasigroups and commutative Moufang loops of exponent 3 (\cite{Fi}; \cite{Ben}, p.133).
In a different, but related direction, if we give a twisted subgroup $T$ of a group $G$
the binary operation $x\star y := xy\iv x$, $x,y\in T$, then $(T,\star)$ is a left
quasigroup which is balanced ($x\star y = y$ iff $y\star x = x$), left distributive
($x\star (y\star z) = (x\star y) \star (x\star z)$), left key ($x\star (x\star y) = y$),
and idempotent ($x\star x = x$). (Other subsets of groups can also be given
this structure, such any conjugacy class with the operation $(x,y) \mapsto xyx\iv$.)
If $T = G$, $(T,\star)$ is called the ``core" of $G$ (this is not the same usage as
in group theory), and the same properties hold even if $G$ is a Moufang loop
\cite{Br}. Studies of these structures, with twisted subgroups as a principal example, can be
found in the work of Nobusawa and his collaborators (see \cite{Nobusawa4} and the
references therein), who were in turn influenced by the work of Loos \cite{Loos} in symmetric
spaces. See also Pierce \cite{Pierce1} \cite{Pierce2} and Umaya \cite{Umaya}. Doro \cite{Doro}
used these structures in his study of simple Moufang loops. Nowadays the structure $(T,\star)$ is
known as an \emph{involutory quandle}, thanks largely to Joyce's applications of the idea to knot
theory \cite{Joyce1}. As far as we have been able to determine, Aschbacher's paper \cite{Asch}
(which was motivated by work of Feder and Vardi \cite{FV}) seems to be the first in which twisted
subgroups are used for a purpose other than the study of quasigroups and loops.
\end{remark}

\section{Bol loops of odd order}
\seclabel{2dbol}

We saw from Example \exref{converse} that a twisted subgroup of
odd order need not be $2$-divisible. However, a twisted subgroup
of odd order which has a compatible Bol loop structure is indeed
$2$-divisible. This is, in fact, a well-known consequence of the
left power-alternative property for Bol loops.

\begin{proposition} (e.g., \cite{Kie})
\proplabel{odd}
The order of any element of a finite Bol loop divides the order of the loop.
\end{proposition}

\noindent In particular, instead of stating results for finite, $2$-divisible Bol loops,
we may simply state them for Bol loops of odd order.

For Bol loops of odd order, the Cauchy and Strong Lagrange theorems
for twisted subgroups immediately transfer to the loop level.

\begin{theorem}
(Cauchy's Theorem)
Let $\LL$ be a Bol loop of odd order. For every prime $p$ dividing $\LL$,
there exists $x\in \LL$ of order $p$.
\end{theorem}

\begin{proof}
By Corollary \cororef{ts-cauchy}, there exists $L(x)\in L(\LL)$ of order
$p$. Since $L(x^n) = L(x)^n$ for all $n$, $x$ has order $p$.
\end{proof}

\begin{remark}
\remlabel{bol-no-hope}
As mentioned in Remark \ref{rem:no-hope} on the twisted subgroup level,
Cauchy's theorem does not extend to all Bol loops, because the simple
Moufang loop of order $180$ does not have an element of order $5$.
\end{remark}

\begin{theorem}
(Strong Lagrange Theorem)
Let $\LL$ be a Bol loop of odd order. If $\KK_1 \subseteq \KK_2
\subseteq \LL$ are subloops, then $|\KK_1|$ divides $|\KK_2|$.
\end{theorem}

\begin{proof}
By Proposition \propref{strongLagrange}, $| L(\KK_1) |$ divides
$| L(\KK_2) |$.
\end{proof}

\begin{problem}
\problabel{strongLagrange}
Does the strong Lagrange property hold for all Bol loops?
\end{problem}

\begin{remark}
If a classification of finite, simple Bol loops were known, then it would be enough
to verify the strong Lagrange property for such loops \cite{CKRV}.
However, this observation merely reduces one hard problem to another.
\end{remark}

\begin{theorem}
Let $\pi$ be a set of odd primes, and let $\LL$ be a finite Bol $\pi$-loop.
Then $\LMlt(\LL)$ is a $\pi$-group.
\end{theorem}

\begin{proof}
This is Theorem \ref{thm:pi} interpreted on the loop level.
\end{proof}

By the Feit-Thompson theorem \cite{FT}, we conclude the following.

\begin{corollary}
If $\LL$ is a finite Bol loop of odd order, then $\LMlt(\LL)$ is a
solvable group.
\end{corollary}

The status of the Sylow and Hall theorems for Bol loops is unclear even
for Bol loops of odd order. For Moufang loops we have the following
results of Glauberman.

\begin{proposition}(Hall's Theorem, \cite{Gl2}, p. 413, Theorem 16 and p. 409, Theorem 12)
Let $\LL$ be a Moufang loop of odd order and let $\pi$ be a set of primes.
Then $\LL$ contains a Hall $\pi$-subloop.
\end{proposition}

The Sylow theorem for Moufang loops of odd order follows immediately,
although Glauberman also gave a separate proof (\cite{Gl2}, p. 410, Theorem 13).
Glauberman also proved other Hall-like properties of $\pi$-subloops
(\cite{Gl2}, p. 413, Theorem 16).

\begin{problem}\hspace*{\fill}
\problabel{bol-Hall}
\begin{enumerate}
\item For a given set $\pi$ of primes, does every Bol loop of odd order have
a Hall $\pi$-subloop?
\item If the answer to (1) is no, then for a given odd prime $p$, does every
Bol loop of odd order have a Sylow $p$-subloop?
\end{enumerate}
\end{problem}

Finally, we turn to some preliminary investigations of simple Bol loops of
odd order. This is motivated by Glauberman's Feit-Thompson theorems
for B-loops and Moufang loops.

\begin{proposition} (\cite{Gl2}, p. 412, Theorem 14 and p. 413, Theorem 16)
\proplabel{solvable}
Let $\LL$ be a Bol loop of odd order. If $\LL$ is a B-loop or
a Moufang loop, then $\LL$ is solvable.
\end{proposition}

\begin{corollary}
\corolabel{extend}
Let $X$ be the class consisting of all Moufang loops of odd order
and all finite B-loops. Let $V$ be any variety of Bol loops of odd order
such that every loop in $V$ is an extension of two loops in $X$. Then
every loop in $V$ is solvable.
\end{corollary}

A left loop is said to have the $A_l$-\emph{property} if every left
inner mapping is an automorphism (\cite{Kie}, p. 35).

\begin{corollary}
\corolabel{Al}
If $\LL$ is an $A_l$ Bol loop of odd order, then $\LL$ is solvable.
\end{corollary}

\begin{proof}
By \cite{FU}, Theorem 4.11, $\LL$ is an extension of a group by a B-loop.
\end{proof}

These considerations pave the way to the following problem. We
will not give a complete answer, but we will present some results
which we think will play a role in its solution.

\begin{problem}
\problabel{solvable}
Do there exist any finite simple Bol loops of odd order?
That is, is every finite Bol loop of odd order solvable?
\end{problem}

Let $\LL$ be a $2$-divisible Bol loop. Since $\LL$ is left power-alternative, that is,
$L(x^n) = L(x)^n$ for all $x\in \LL$, we may use the $2$-divisible twisted subgroup
$L(\LL)$ to define a B-loop operation on $\LL$. For $x,y\in \LL$, we have
$L(x)\odot L(y) = (L(x)L(y)^2L(x))^{1/2} = L((x\cdot (y\cdot x))^{1/2})$. This us
leads us to the following.

\begin{definition}
Let $\LL$ be a $2$-divisible Bol loop. The \emph{B-loop associated to} $\LL$ is
$(\LL,\odot)$ with the binary operation $\odot : \LL\times \LL \to \LL$ given by
$x\odot y = (x\cdot ((y\cdot y)\cdot x))^{1/2}$. We will denote the B-loop $(\LL,\odot)$
by $\LL(1/2)$, and we follow a similar convention for subloops. Left multiplication
maps for $\LL(1/2)$ will be denoted by $M(x)y := x\odot y$.
\end{definition}

\begin{remark}
In case $\LL$ is already a B-loop, $(\LL,\odot)$ is just $\LL$ itself. This is because
every Bruck loop satisfies the identity $(x\cdot y)^2 = x \cdot (y^2 \cdot x)$
(\cite{Kie}, p. 73, 6.8(1)).
\end{remark}

The B-loop associated to a Moufang loop of odd order was the key component
in Glauberman's proofs of the Hall, Sylow, and Feit-Thompson theorems in
\cite{Gl2}. The idea was to ``pull back" the results from the associated B-loop
to the Moufang loop. Since arbitrary Bol loops are not as structured as Moufang
loops, one cannot expect this idea to work quite so well. Nevertheless, we can make
some progress.

\begin{lemma}
\lemlabel{PisB}
Let $\LL$ be a $2$-divisible Bol loop. Then the squaring map
$s:\LL \to \LL; x\mapsto x\cdot x$ conjugates $\LMlt(\LL(1/2))$
to $\PMlt(\LL)$ in $\LL !$.
\end{lemma}

\begin{proof}
For each $x\in \LL$, $M(x) = s\iv P(x) s$.
\end{proof}

\begin{theorem}
\thmlabel{rep2}
Let $\LL$ be a $2$-divisible Bol loop, and let $G = \LMlt(\LL)$.
The mapping $L(\LL)\to M(\LL(1/2)); L(x)\mapsto M(x)$ extends
to a homomorphism from $G$ onto $\LMlt(\LL(1/2))$ if and only
if $\LL^{\prime} = \{ 1\}$. The kernel of the homomorphism is
$Z(G)\cap C_G(\tau)$ where $\tau\in\Aut(G)$ is the Aschbacher
automorphism. 
\end{theorem}

\begin{proof}
This follows from Lemma \lemref{PisB} and Theorem \thmref{rep1}.
\end{proof}

\begin{lemma}
\lemlabel{normal}
A subloop $\KK$ of a $2$-divisible Bol loop $\LL$ is normal
if and only if, for all $x,y\in \LL$,
$x \cdot (\KK\cdot y) = \KK \cdot (x\cdot y)$.
\end{lemma}

\begin{proof}
Referring to the conditions in \peqref{normality},
we see that only one direction requires proof. Thus assume
$x \cdot (\KK\cdot y) = \KK \cdot (x\cdot y)$ for
all $x,y \in \LL$. Fix $x\in \LL$ and set $u = x^{1/2}$.
Using LAP and the left Bol identity,
\[
\begin{array}{rcl}
x \cdot (\KK \cdot y) &=& u\cdot (u\cdot (\KK\cdot y)) = u\cdot (\KK\cdot (u\cdot y))\\
&=& (u\cdot (\KK\cdot u))\cdot y = (u\cdot (u\cdot \KK))\cdot y = (x\cdot \KK)\cdot y.
\end{array}
\]
Thus the other condition of \peqref{normality} holds, and so $\KK$ is normal.
\end{proof}

Our next result is inspired by Aschbacher's normality condition
for subloops (\cite{Asch2}, condition (NC)). It enables us to express
normality directly in terms of the left multiplication group.

\begin{lemma}
\lemlabel{normal2}
A subloop $\KK$ of a $2$-divisible Bol loop $\LL$ is normal if and only if,
for each $x\in \LL$, $y\in \KK$, $g\in \LMlt(\LL)$,
\begin{equation}
\eqlabel{normal}
L(x) L(y) L(x\iv) = L(z) ghg\iv
\end{equation}
for some $z\in \KK$, $h\in \LMlt_1(\LL)$.
\end{lemma}

\begin{proof}
Set $G = \LMlt(\LL)$, $H = \LMlt_1(\LL)$. Fix $x\in \LL$, $y\in \KK$.
Since $L(\LL)$ is a transversal of each conjugate $g H g\iv$, $g\in G$,
we have $L(x) L(y) L(x\iv) = L(z) ghg\iv$ for some $z\in \LL$, $h\in H$. Applying
both sides to $w = g1$, we have $x\cdot (y\cdot (x\iv \cdot w)) = z\cdot w$.
Now if $\KK$ is normal, then by Lemma \lemref{normal} (or just \peqref{normality}),
$x\cdot (y\cdot (x\iv \cdot w)) = u\cdot w$ for some $u\in \KK$. Thus
$z = u$ and so $z \in \KK$ so that \peqref{normal} holds. Conversely, if \peqref{normal} holds,
then fix $v\in \LL$ and set $g = L(v)$. Let $z\in \KK$, which depends on $g$, be given
as in \peqref{normal}. Apply both sides of \peqref{normal} to $z$ to get
$x\cdot (y\cdot (x\iv \cdot v)) = z\cdot v$. By Lemma \lemref{normal}, 
$\KK$ is normal.
\end{proof}

In the Moufang case, the following result is (\cite{Gl2}, p. 401, Lemma 7(b)).
The proof in the general case is essentially the same, but with care in
the parenthesization.

\begin{lemma}
\lemlabel{bol-subloops}
Let $\LL$ be a $2$-divisible Bol loop and suppose $\KK$ is a
subloop of $\LL(1/2)$. Then $\KK$ is a subloop of $\LL$ if
and only if $x\iv\cdot (\KK\cdot x) = \KK$ for all $x\in \KK$.
\end{lemma}

\begin{proof}
The ``only if" is obvious, so assume $x\iv\cdot (\KK\cdot x) = \KK$
for all $x\in \KK$. If $x\in \KK$, then $\langle x\rangle \subset \KK$, 
because powers in $\LL$ agree with powers in $\LL(1/2)$. Now fix $x,y\in \KK$
and set $u = x^{1/2}$, $v = y^{1/2}$. Then using the definition of $\odot$, the Bol
identity, and LAP, $\KK$ contains
$u\cdot ((u\odot v)^2 \cdot u\iv) = u\cdot ((u\cdot (v^2 \cdot u))\cdot u\iv)
= u^2 \cdot v^2 = x\cdot y$. A subset of a Bol loop closed under inversion
and multiplication is a subloop (see, e.g., \cite{Kie}, p. 50, 3.10(4)).
\end{proof}

\begin{lemma}
\lemlabel{main}
Let $\LL$ be a radical-free, $2$-divisible Bol loop, let $G = \LMlt(\LL)$,
let $\tau\in\Aut(G)$ be the Aschbacher automorphism, and assume that
$Z(G)\cap C_G(\tau) = \langle 1\rangle$. If $\KK(1/2)$ is a normal subloop of
$\LL(1/2)$, then $\KK$ is a normal subloop of $\LL$.
\end{lemma}

\begin{proof}
Fix $x\in \LL$, $y\in \KK$, $g\in G$. As in the proof of
Lemma \lemref{normal2}, there exists $z\in \LL$, $h\in \LMlt_1(\LL)$
such that
\begin{align*}
L(x) L(y) L(x\iv) = L(z) ghg\iv \tag{*}
\end{align*} 
The hypotheses and Theorem \ref{thm:rep2} imply 
$\LMlt(\LL(1/2)) \cong G$, and that on generators, the isomorphism
sends each $L(x)$ to $M(x)$. Applying this to (*), we have
$M(x) M(y) M(x\iv) = M(z) \tilde{g}\tilde{h}\tilde{g}\iv$ for
some $\tilde{g}\in \LMlt(\LL(1/2))$, $\tilde{h}\in \LMlt_1(\LL(1/2))$.
If $\KK(1/2)$ is normal, then by Lemma \lemref{normal2}, $z\in \KK$.
Thus by (*), the condition of Lemma \lemref{normal2} is satisfied for
the sub\emph{set} $\KK$ of $\LL$, and all that remains is to show
that $\KK$ is a subloop. Taking $x\in \KK$ and $g = 1$ in (*), and
applying both sides to $1$, we have that for each $x,y\in \KK$, there
exists $z\in \KK$ such that $x\cdot (y\cdot x\iv) = z$. 
Now apply Lemma \lemref{bol-subloops}.
\end{proof}

\begin{theorem}
\thmlabel{main}
Let $\LL$ be a simple Bol loop of odd order, let $G = \LMlt(\LL)$,
and (since $\LL$ is radical-free), let $\tau\in\Aut(G)$
denote the Aschbacher automorphism. Then
$Z(G)\cap C_G(\tau) \neq \langle 1 \rangle$.
\end{theorem}

\begin{proof}
By Proposition \propref{solvable}, $\LL(1/2)$ has a nontrivial normal
subloop $\KK(1/2)$. If $Z(G)\cap C_G(\tau) = \langle 1 \rangle$, then
Lemma \lemref{main} implies that $\KK$ is a nontrivial normal subloop
of $\LL$.
\end{proof}

\begin{corollary}
\corolabel{main}
Let $\LL$ be a simple Bol loop of odd order, and let $G = \LMlt(\LL)$.
Then $G$ has nontrivial center, and $\Nuc_r(\LL)$ contains an abelian
subgroup $\MM = \{ x\in \Nuc_r(\LL) : R(x)\in G\} \neq \langle 1\rangle$.
\end{corollary}

\begin{proof}
This follows from Theorem \thmref{main} and Lemma \lemref{kernel}.
\end{proof}

\begin{remark}
\remlabel{mfg}
As a corollary, we obtain a new proof of the Moufang part of Proposition
\propref{solvable}. In a Moufang loop, the three nuclei agree and the
nucleus is normal (\cite{Pf}, p. 90, Corollary IV.1.5). Thus a simple Moufang loop
has trivial nucleus, and so by Corollary \cororef{main}, cannot have odd order.
\end{remark}

In general, the right nucleus of a left Bol loop need not be a normal subloop.
This follows from a construction of D.~Robinson and K.~Robinson \cite{RR},
translated from right Bol loops to left Bol loops. However, their construction gives
a Bol loop of even order. Thus the following problem still seems to be open,
even for B-loops.

\begin{problem}
Does there exist a Bol loop of odd order with nonnormal right nucleus?
\end{problem}

\begin{acknowledgment}
We thank Michael~Aschbacher, James~Beidleman, George~Glauberman, and
Derek~Robinson for helpful comments.
\end{acknowledgment}


\begin{thebibliography}{99}
\bibitem{Asch}
M.~Aschbacher,
Near subgroups of finite groups,
\textit{J. Group Theory} \textbf{1} (1998) 113-129.

\bibitem{Asch2}
M.~Aschbacher,
Bol loops of exponent 2, 
to appear in \textit{J. Algebra}.

\bibitem{Baer}
R.~Baer,
Nets and groups,
\textit{Trans. Amer. Math. Soc.} \textbf{47} (1939), 110-141.

\bibitem{Bel1}
V.D.~Belousov,
The core of a Bol loop,
in \textit{Collection in General Algebra (Sem)}, 53-66, Akad. Nauk.
Moldav. SSR, Kishinev, 1965 (Russian)

\bibitem{Bel2}
V.D.~Belousov,
\textit{Foundations of the Theory of Quasigroups and Loops},
Izdat. Nauka, Moscow, 1967 (Russian)

\bibitem{Ben}
L.~B\'{e}n\'{e}teau,
Commutative Moufang loops and related groupoids,
Chapter IV in \cite{CPS}, 115-142

\bibitem{Br}
R.H.~Bruck,
\textit{A Survey of Binary Systems}
Springer Verlag, Berlin, 1971

\bibitem{CPS}
O.~Chein, H.O.~Pflugfelder, and J.D.H.~Smith (eds.),
\textit{Quasigroups and Loops: Theory and Applications},
Sigma Series in Pure Math. \textbf {9}, Heldermann Verlag, Berlin, 1990

\bibitem{CKRV}
O.~Chein, M.K.~Kinyon, A.~Rajah, and P.~Vojt\v{e}chovsk\'{y},
Loops and the Lagrange property,
\textit{Results in Math.} \textbf{43} (2003), 74-78.

\bibitem{Doro}
S.~Doro,
Simple Moufang loops,
\textit{Math. Proc. Cambridge} \textbf{83} (1978), 377-392.

\bibitem{FV}
T.~Feder and M.~Vardi,
The computational structure of monotone monadic SNP and constraint
satisfaction: a study through Datalog and group theory,
\textit{SIAM J. Computing} \textbf{28} (1998) 57-104.

\bibitem{Feder}
T.~Feder,
Strong near subgroups and left gyrogroups,
\textit{J. Algebra}  \textbf{259}  (2003)  177-190

\bibitem{FT}
W.~Feit and J.~Thompson,
Solvability of groups of odd order,
\textit{Pacific J. Math.} \textbf{13} (1963), 775-1029.

\bibitem{Fi}
B.~Fischer,
Distributive Quasigruppen endlicher Ordnung,
\textit{Math. Zeit.} \textbf{83} (1964) 267-303.

\bibitem{FU}
T.~Foguel and A.A.~Ungar,
Involutory decomposition of groups into twisted subgroups and subgroups,
\textit{J. Group Theory} \textbf{3} (2000) 27-46.

\bibitem{Foguel1}
T.~Foguel,
Groups with polar decompositions,
\textit{Results in Math.}  \textbf{42}  (2002)  69-73

\bibitem{FN}
M.~Funk and P.T.~Nagy,
On collineation groups generated by Bol reflections,
\textit{J. Geometry} \textbf{48} (1993) 63-78.

\bibitem{Gl1}
G.~Glauberman,
On loops of odd order I,
\textit{J. Algebra} \textbf{1} (1964) 374-396.

\bibitem{Gl2}
G.~Glauberman,
On loops of odd order II,
\textit{J. Algebra} \textbf{8} (1968) 393-414.

\bibitem{Gl3}
G.~Glauberman,
Central elements in core-free groups,
\textit{J. Algebra} \textbf{4} (1966) 403-420.

\bibitem{Joyce1}
D.~Joyce,
A classifying invariant of knots, the knot quandle,
\textit{J. Pure Appl. Alg.} \textbf{23} (1982), 37-66.

\bibitem{Kie}
H.~Kiechle,
\textit{Theory of K-loops},
Lecture Notes in Mathematics 1778, Springer-Verlag,
Berlin Heidelberg New York, 2002.

\bibitem{Kik}
M.~Kikkawa,
On some quasigroups of algebraic models of symmetric spaces II,
\textit{Mem. Fac. Lit. Sci., Shimane Univ., Nat. Sci.} \textbf{7} (1974) 29-35.

\bibitem{KJ}
M.K.~Kinyon and O.~Jones,
Loops and semidirect products,
\textit{Comm. Algebra} \textbf{28} (2000) 4137-4164.

\bibitem{Kr}
A.~Kreuzer,
Beispiele endlicher und unendlicher K-loops,
\textit{Results in Math.} \textbf{23} (1993) 355-362.

\bibitem{Loos}
O.~Loos,
\textit{Symmetric spaces I.}
J. Benjamin, New York, 1969.

\bibitem{Na2}
G.P.~Nagy,
Group invariants of certain Burn loop classes,
\textit{Bull. Belg. Math. Soc.} \textbf{5} (1998) 403-415.

\bibitem{Na3}
G.P.~Nagy,
\textit{Bol-t\"{u}kr\"{o}z\'{e}sek alkalmaz\'{a}sai},
Ph.D. dissertation, Bolyai Institute, Szeged, 1999 (Hungarian).

\bibitem{NS}
P.T.~Nagy and K.~Strambach,
Loops, their cores and symmetric spaces,
\textit{Israel J. Math.} \textbf{105} (1998) 285-322.

\bibitem{Nobusawa4}
N.~Nobusawa,
Orthogonal groups and symmetric sets,
\textit{Osaka J. Math.} \textbf{20} (1983), 5-8.

\bibitem{Pf}
H.O.~Pflugfelder,
\textit{Quasigroups and Loops: An Introduction},
Sigma Series in Pure Math. \textbf{8},  Heldermann Verlag, Berlin, 1990.

\bibitem{Ph}
J.D.~Phillips,
Quotients of groups,
\textit{Tamkang Journal of Mathematics}
\textbf{28} (1997) 271-275.

\bibitem{Pierce1}
R.~S.~Pierce,
Symmetric groupoids,
\textit{Osaka J. Math.} \textbf{15} (1978), 51-76.

\bibitem{Pierce2}
R.~S.~Pierce,
Symmetric groupoids II,
\textit{Osaka J. Math.} \textbf{16} (1979), 317-348.

\bibitem{OpenProblems}
\textit{Problems in Loop Theory and Quasigroup Theory},
at \\
{\tt http://adela.karlin.mff.cuni.cz/\symbol{126}loops03/plq/main.html}


\bibitem{Ro}
D.A.~Robinson,
Bol loops,
\textit{Trans. Amer. Math. Soc.} \textbf{123} (1966), 341-354.

\bibitem{RR}
D.A.~Robinson and K.~Robinson,
\textit{Arch. Math. (Basel)} \textbf{61} (1993), 596-600.
 
\bibitem{Umaya}
N.~Umaya,
On symmetric structure of a group.
\textit{Proc. Japan Acad.} \textbf{52} (1976), 174-176.

\bibitem{Ung}
A.A.~Ungar,
Thomas precession: its underlying gyrogroup
axioms and their use in hyperbolic geometry
and relativistic physics,
\textit{Found. Phys.} \textbf{27} (1997), 881-951.

\end{thebibliography}
\end{document}